\documentclass[11pt]{article}

\usepackage{amscd,amsmath, amssymb, fancyhdr}

\newcommand{\version}{version 1.01,\ \   June 12, 2009}

\newcommand{\la}{\lambda}
\newcommand{\al}{\alpha}
\newcommand{\be}{\beta}

\newcommand{\tom}{\tilde \omega}

\newcommand{\Ll}{\operatorname{Lie}}

\newcommand{\RR}{\mathbb{R}}
\newcommand{\ZZ}{\mathbb{Z}}

\numberwithin{equation}{section}

\def\eqref#1{(\ref{#1})}

\newcommand{\arrow}{{\:\longrightarrow\:}}
\newcommand{\Z}{{\Bbb Z}}
\newcommand{\C}{{\Bbb C}}
\newcommand{\R}{{\Bbb R}}

\def\1{\sqrt{-1}\:}
\newcommand{\restrict}[1]{{\left|_{{\phantom{|}\!\!}_{#1}}\right.}}
\newcommand{\cntrct}                
{\hspace{2pt}\raisebox{1pt}{\text{$\lrcorner$}}\hspace{2pt}}

\renewcommand{\phi}{\varphi}
\renewcommand{\epsilon}{\varepsilon}
\renewcommand{\geq}{\geqslant}

\newcounter{Mycounter}[section]
\newcounter{lemma}[section]
\newcounter{claim}[section]
\newcounter{sublemma}[section]
\newcounter{corollary}[section]
\renewcommand{\thecorollary}{{Corollary \thesection.\arabic{corollary}}}
\newcommand{\corollary}{%
     \setcounter{corollary}{\value{Mycounter}}
     \refstepcounter{corollary}
     \stepcounter{Mycounter}
     {\noindent \bf \thecorollary:\ }}

\newcounter{theorem}[section]
\renewcommand{\thetheorem}{{Theorem \thesection.\arabic{theorem}}}
\newcommand{\theorem}{%
     \setcounter{theorem}{\value{Mycounter}}
     \refstepcounter{theorem}
     \stepcounter{Mycounter}
     {\noindent \bf \thetheorem:\ }}

\newcounter{conjecture}[section]
\newcounter{proposition}[section]
\renewcommand{\theproposition}
       {{Proposition \thesection.\arabic{proposition}}}
\newcommand{\proposition}{%
     \setcounter{proposition}{\value{Mycounter}}
     \refstepcounter{proposition}
     \stepcounter{Mycounter}
     {\noindent \bf \theproposition:\ }}

\newcounter{definition}[section]
\renewcommand{\thedefinition}
       {{Definition~\thesection.\arabic{definition}}}
\newcommand{\definition}{%
     \setcounter{definition}{\value{Mycounter}}
     \refstepcounter{definition}
     \stepcounter{Mycounter}
     {\noindent \bf \thedefinition:\ }}

\newcounter{example}[section]
\newcounter{remark}[section]
\renewcommand{\theremark}{{Remark \thesection.\arabic{remark}}}
\newcommand{\remark}{%
     \setcounter{remark}{\value{Mycounter}}
     \refstepcounter{remark}
     \stepcounter{Mycounter}
     {\noindent \bf \theremark:\ }}

\newcounter{problem}[section]
\newcounter{question}[section]
\makeatletter

\@addtoreset{equation}{section}
\@addtoreset{footnote}{section}
\makeatother

\def\blacksquare{\hbox{\vrule width 5pt height 5pt depth 0pt}}
\def\endproof{\hfill\blacksquare}

\begin{document}
\begin{center}
{\LARGE\bf
Automorphisms of locally conformally\\[3mm] K\"ahler manifolds}\\[3mm]

Liviu Ornea\footnote{Partially supported by a PN2-IDEI grant, nr. 529.}
and Misha Verbitsky\footnote{Partially supported by the
grant RFBR for support of scientific schools
NSh-3036.2008.2 and RFBR grant 09-01-00242-a.\\


}

\end{center}

{\small
\hspace{0.15\linewidth}
\begin{minipage}[t]{0.7\linewidth}
{\bf Abstract} \\
A manifold $M$ is locally conformally K\"ahler (LCK)
if it admits a  K\"ahler covering $\tilde M$ with
monodromy acting by holomorphic homotheties. For
a compact connected group $G$ acting on an LCK
manifold by holomorphic automorphisms, an averaging
procedure gives a $G$-invariant LCK metric.
Suppose that $S^1$ acts on an LCK manifold $M$ by
holomorphic isometries, and the lifting of this action
to the K\"ahler cover $\tilde M$ is not isometric.
We show that $\tilde M$ admits an automorphic
K\"ahler potential, and hence (for
$\dim_\C M >2$) the manifold $M$
can be embedded to a Hopf manifold.
\end{minipage}
}

\tableofcontents

\section{Introduction}

\subsection{Locally conformally K{\"a}hler manifolds}

Locally conformally K{\"a}hler (LCK)
manifolds are, by definition,
complex manifolds of $\dim_\C>1$
admitting a K{\"a}hler covering
with deck transformations acting by K{\"a}hler homotheties. We shall
usually denote with $\tilde\omega$ the K\"ahler form on the covering.

An equivalent definition, at the level of the manifold itself, postulates
the existence of an open covering $\{U_\al\}$ with  local K\"ahler
metrics $g_\al$. It requires that on overlaps $U_\al\cap U_\be$, these
local K\"ahler metrics are homothetic: $g_\al=c_{\al\be}g_\be$.
The metrics
$e^{f_\al}g_\al$ glue to a global metric whose associated two-form
$\omega$ satisfies the integrability condition
$d\omega=\theta\wedge\omega$, thus being locally conformal with the
K\"ahler metrics $g_\al$. Here  $\theta\restrict{U_\al}=df_\al$. The
closed 1-form
$\theta$, which represents the cocycle $c_{\al\be}$, is called {\bf the
Lee form}. Obviously, any other representative
of this cocycle, $\theta'=\theta+dh$, produces another LCK metric,
conformal with the initial one. This gives
another definition of an LCK structure, which
will be used in this paper.

\hfill

\definition
Let $(M, \omega)$ be a complex Hermitian manifold,
$\dim_\C M >1$, with $d\omega = \theta\wedge\omega$,
where $\theta$ is a closed 1-form. Then
$M$ is called {\bf a locally conformally K{\"a}hler} (LCK)
manifold

\hfill

We refer to \cite{drag} for an overview and to
\cite{_OV_Top_Potential_} for more recent results.

\subsection{Bott-Chern cohomology and automorphic potential}

Let $(\tilde M, \tilde\omega)$ be a K\"ahler covering
of an LCK manifold $M$, and let $\Gamma$ be the deck transform
group of $[\tilde M:M]$.
Denote by $\chi:\; \Gamma \arrow \R^{>0}$
the corresponding character of $\Gamma$,
defined through the scale factor of $\tilde \omega$:
\[
  \gamma^*\tilde\omega=\chi(\gamma) \tilde\omega, \ \
  \forall \gamma\in\Gamma.
\]

\definition
A differential form $\alpha$ on $\tilde M$ is called
{\bf automorphic} if $\gamma^*\alpha =
\chi(\gamma)\alpha$, where  $\chi:\; \Gamma \arrow \R^{>0}$
is the character of $\Gamma$ defined above.

\hfill

A useful tool in the study of LCK geometry is the weight bundle $L\arrow
M$. It is a topologically trivial line bundle, associated
to the representation $\mathrm{GL}(2n,\mathbb{R})\ni A \mapsto | \det
A|^{\frac{1}{n}}$, with flat connection
defined as $D:= \nabla_0 + \theta$, where $\nabla_0$
is the trivial connection. It allows regarding automorphic objects on
$\tilde M$ as objects on $M$ with values in $L$.

\hfill

\definition
Let $M$ be an LCK manifold, $\Lambda^{1,1}_{\chi,
  d}(\tilde M)$ the space of closed, automorphic
$(1,1)$-forms on its K\"ahler covering $\tilde M$, and
let $\mathcal{C}^\infty_\chi(\tilde M)$ be
the space of automorphic functions on $\tilde M$.
Consider the quotient
\[
H^{1,1}_{BC}(M, L):= \frac{\Lambda^{1,1}_{\chi,
  d}(\tilde M)}{dd^c(\mathcal{C}^\infty_\chi(\tilde M))},
\]
where $d^c=-IdI$.
This group is finite-dimensional. It is
called {\bf the Bott-Chern cohomology group of an LCK manifold}
(for more details, see  \cite{_ov:MN_}). It
is independent from the choice of
the covering $\tilde M$.

\hfill

\remark
The K\"ahler form $\tilde\omega$ on $\tilde M$
is obviously closed and automorphic. Its cohomology
class $[\tilde \omega]\in H^{1,1}_{BC}(M, L)$
is called {\bf the Bott-Chern class of $M$}.
It is an important cohomology invariant of
an LCK manifold, which can be considered as
an LCK analogue of the K\"ahler class.

\hfill

\definition
Let $(\tilde M, \tilde\omega)$ be a K\"ahler covering
of an LCK manifold $M$. We say that $M$ is
{\bf an LCK manifold with an automorphic potential}
if $\tilde \omega = dd^c \phi$, for some automorphic
function $\phi$ on $\tilde M$. Equivalently, $M$ is
an LCK manifold with an automorphic potential,
if its Bott-Chern class vanishes.

\hfill

Compact LCK manifolds with automorphic potential are embeddable in
Hopf manifolds, see \cite{_ov:MN_}. The existence
of an automorphic potential leads to important topological
restrictions on the fundamental group, see
\cite{_OV_Top_Potential_} and \cite{_Kokarev_Kotschick:Fibrations_}.

The class of compact complex manifolds admitting
an LCK metric with automorphic potential is
stable under small complex deformation, \cite{_OV:Potential_}.
This statement should be considered as an LCK
analogue of Kodaira's celebrated K\"ahler
stability theorem.
The only way (known to us) to construct LCK metrics
on some non-Vaisman manifolds, such as the
Hopf manifolds not admitting a Vaisman structure,
is by deformation, applying the stability of
automorphic potential under small deformations.

\subsection{Automorphisms of LCK and Vaisman manifolds}

\definition
A {\bf Vaisman manifold} is an LCK manifold $(M, \omega, \theta)$
with $\nabla \theta=0$, where $\theta$ is its Lee form,
and $\nabla$ the Levi-Civita connection.

\hfill

As shown e.g. in \cite{_Verbitsky_vanishing_}, a
Vaisman manifold has an automorphic potential,
which can be written down explicitely as
$\tilde \omega(\pi^* \theta, \pi^*\theta)$,
where $\pi^*\theta$ is the lift of the Lee
form to the considered K\"ahler covering of $M$.

\hfill

Compact Vaisman manifolds can be characterized in terms of
their automorphisms group.

\hfill

\theorem\label{_Kami_Or_Theorem_}
(\cite{_Kamishima_Ornea_})
Let $(M, \omega)$ be a compact LCK manifold admitting a
holomorphic, conformal action of $\C$ which
lifts to an action by non-trivial homotheties on its
K\"ahler covering. Then $(M, \omega)$ is conformally
equivalent to a Vaisman manifold. \endproof

\hfill

 Other properties of the various
transformations groups of LCK manifolds were studied in
\cite{_moroianu_ornea_} and \cite{_GOP:reduction_}.

It was proven in \cite{_OV_Top_Potential_}
that any compact LCK manifold
with automorphic  potential can be obtained as a
deformation of a Vaisman manifold.  Many of the known examples of LCK
manifolds
are Vaisman (see \cite{belgun} for a complete list of Vaisman compact
complex surfaces), but there are also non-Vaisman ones: one of the Inoue
surfaces (see \cite{belgun}, \cite{_Tricerri_}), its higher-dimensional
generalization in \cite{ot}, and the new examples found in \cite{FP} on
parabolic and hyperbolic Inoue surfaces. Also, a blow-up
of a Vaisman manifold is still LCK (see \cite{_Tricerri_}, \cite{vuli}),
but not Vaisman,
and has no automorphic potential.

In this paper, we show that
LCK manifolds with automorphic
potential can be characterized in terms of existence of a
particular subgroup of automorphisms. In Section 2,
we prove the following theorem.

\hfill

\theorem\label{_S^1_main_Theorem_}
Let $M$ be a compact complex manifold, equipped
with a holomorphic $S^1$-action and an LCK metric
(not necessarily compatible). Suppose that
the weight bundle $L$, restricted to a general
orbit of this $S^1$-action, is non-trivial
as a 1-dimensional local system. Then $M$ admits
an LCK metric with an automorphic potential.

\hfill

\remark
The converse statement seems to be true as well. We conjecture that
given a LCK manifold $M$ with a automorphic
potential, $M$ always admits a holomorphic $S^1$-action
of this kind. To motivate this conjecture, consider a Hopf
manifold $M$ (Hopf manifolds are known to admit an LCK metric
with an automorphic potential, see e.g. \cite{_OV:Potential_}).
Suppose that $M$ is a quotient of
$\C^n\backslash 0$ by a group $\Z$ acting by linear
contractions, $M=\C^n\backslash 0/\langle A\rangle$,
with $A$ a linear operator with all eigenvalues $\alpha_i$
satisfying $|\alpha_i|<1$.\footnote{Such Hopf manifolds are
called {\bf linear}.} Then the
holomorphic diffeomorphism
flow associated with the vector field $\log A$ leads
to a holomorphic $S^1$-action on $M$.

\hfill

\remark
\ref{_Kami_Or_Theorem_}
implies that an LCK manifold $M$ with a certain conformal action
of $\C$ is conformally equivalent to a Vaisman manifold.
By contrast, \ref{_S^1_main_Theorem_} does not postulate that
the given $S^1$-action is compatible with the metric.
Neither does  \ref{_S^1_main_Theorem_} say anything about
the given LCK metric on $M$. Instead, \ref{_S^1_main_Theorem_}
says that some other LCK structure on the same complex
manifold has an automorphic potential. This new metric
is obtained (see Subsection \ref{_avera_second_Subsection_})
by a kind of convolution, by averaging the old one with
some weight function, wich depends on the cohomological
nature of the $S^1$-action. In particular, the original
LCK metric may have no potential. In \cite[Conjecture 6.3]{_ov:MN_}
it was conjectured that all LCK metrics on a Vaisman
manifold have potential; this conjecture is still unsolved.

\hfill

As shown in \cite{_ov:MN_} and \cite{_OV_Top_Potential_},
\ref{_S^1_main_Theorem_} implies the following corollary.

\hfill

\corollary
Let $M$ be a compact LCK manifold  of complex dimension
$n\geq 3$. Suppose that
the weight bundle $L$ restricted to a general
orbit of this $S^1$-action is non-trivial
as a 1-dimensional local system. Then $\tilde M$
is diffeomorphic to a Vaisman manifold, and admits
a holomorphic embedding to a Hopf manifold.


\section{The proof of the main theorem}


\subsection{Averaging on a compact transformation group}

For the sake of
completeness, we recall the following procedure described in the proof of
\cite[Th. 6.1]{_ov:MN_}. Let $G$ be a compact subgroup of
$\mathrm{Aut}(M)$. Averaging the Lee form $\theta$ on $G$, we obtain a
closed $1$-form $\theta'$ which is $G$-invariant and stays in the same
cohomology class as $\theta$: $\theta'=\theta+df$. Then
$\omega'=e^{-f}\omega$ is a LCK form with Lee form $\theta'$ and conformal
to $\omega$. Hence, we may  assume from the beginning that
$\theta$ (corresponding to $\omega$) is $G$-invariant.

Now, for any $a\in G$, $a^*\omega$ satisfies
\begin{equation}\label{_rho^*_theta_inva_Equation_}
   d(a^*\omega) = a^*\omega\wedge a^*\theta
   = a^*\omega\wedge \theta.
\end{equation}
Averaging  $\omega$ over $G$ and applying
\eqref{_rho^*_theta_inva_Equation_}, we find a
$G$-invariant Hermitian form $\omega'$ which satisfies
\[
d\omega' = \omega'\wedge \theta.
\]
Therefore, we may also assume that $\omega$ is $G$-invariant.

In conclusion, by averaging on $S^1$, we obtain a new LCK metric,
conformal with the
initial one, w.r.t. which $S^1$ acts by (holomorphic) isometries and whose
Lee form is $S^1$-invariant. Hence,
we may suppose from the beginning that $S^1$ acts by holomorphic
isometries of the given LCK metric.

This implies that the lifted action of $\RR$ acts by homotheties of
the global K\"ahler metric with K\"ahler form $\tilde \omega  $.
Indeed, $a^*\tilde \omega= f\tilde \omega$, but
$d(a^*\tilde\omega)=0=df\wedge\tilde\omega$,
and  multiplication by $\tilde \omega$
is injective on $\Lambda^1(M)$, as $\dim_\C M>1$,
hence $df\wedge\tilde\omega=0$ implies $df=0$.

The monodromy of the weight bundle along
an orbit $S$ of the $S^1$-action can be computed
as $\int_S \theta$, hence this monodromy is not changed
by the averaging procedure. Therefore,
it suffices to prove \ref{_S^1_main_Theorem_}
assuming that $\omega$ is $S^1$-invariant.

In this case, the lift of the $S^1$-action
on $\tilde M$ acts on the K\"ahler form
$\tilde M$ by homotheties, and the corresponding
conformal constant is equal to the
monodromy of $L$ along the orbits of
$S^1$. Therefore, we may assume that
$S^1$ is lifted to an $\R$ action on
$\tilde M$ by non-trivial homotheties.

\subsection{The main formula}

Let now $A$ be the vector field on
$\tilde M$ generated by the
$\RR$-action. Then $A$ is holomorphic and homothetic, i.e.
$$\Ll_A\tilde \omega  =\la\tilde \omega  , \quad \la\in\RR^{>0}.$$
Denote:
$$A^c=IA,\quad \eta=A\cntrct  \tilde \omega  , \quad \eta^c=I\eta.$$
Note that, by definition,
$(I\al)(X_1,\ldots,X_k)=(-1)^k\al(IX_1,\ldots,IX_k).$

We now
prove the following formula, which is the key
to the rest of our argument.

\hfill

\proposition
Let $A$ be a vector field acting on a
K\"ahler manifold $\tilde M$ by
holomorphic homotheties: $\Ll_A\tilde \omega  =\la\tilde\omega$.
Then
\begin{equation}\label{unu}
 dd^c|A|^2= \la^2\tilde \omega  +\Ll_{A^c}^2\tilde \omega,
\end{equation}
where $A^c = I(A)$.

\hfill

{\bf Proof:} Replacing $A$ by $\lambda^{-1}A$, we may
assume that $\lambda=1$. By Cartan's formula,
$$\Ll_A\tilde \omega  =d(A\cntrct \tom)=d\eta,$$
and hence, as $\eta(A)=0$,
$$\Ll_A\eta=d(A\cntrct \eta)+A\cntrct  d\eta=A\cntrct (\tilde \omega
)=\eta.$$
As $A$ is holomorphic, this implies $\Ll_A\eta^c=\eta^c$. But, again
with Cartan's formula:
$$\Ll_A\eta^c=d(A\cntrct  \eta^c)+A\cntrct  d\eta^c=-d|A|^2+A\cntrct
d\eta^c.$$
Hence:
$$d^cd|A|^2=- d^c\eta^c+d^c(A\cntrct  d\eta^c),$$
We note that:
$$d^c\eta^c=-Id\eta= -I\tilde \omega  =\tilde \omega  ,$$
as $\tilde \omega  $ is $(1,1)$. Then, to compute $d^c(A\cntrct
d\eta^c)$, observe first that
$$\Ll_{A^c}\tilde \omega  =d(IA\cntrct \tilde \omega  )=d\tilde \omega
(IA,\cdot)=d\eta^c.$$
Thus, as $\tilde \omega  $ and $\Ll_{A^c}\tilde \omega  $ are $(1,1)$,
and by Cartan's formula again:
\begin{align*}
d^c(A\cntrct  d\eta^c)&=-IdI(A\cntrct  \Ll_{A^c}\tilde \omega) =Id(A^c\cntrct
\Ll_{A^c}\tilde \omega  )\\
&=I\Ll_{A^c}^2\tilde \omega =-\Ll_{A^c}^2\tilde \omega .
\end{align*}
This proves \eqref{unu}.

\subsection{The second averaging argument}
\label{_avera_second_Subsection_}

Clearly, the action of the Lie derivative on
$\Omega^\bullet(M)$ can be
extended to the Bott-Chern cohomology groups by $\Ll_X[\al]=[\Ll_X\al]$. Then
\eqref{unu} tells us that
$$\Ll_{A^c}^2[\tilde\omega]=-\la^2[\tilde\omega],$$
where $[\tilde\omega]$ is the class of $\tom$ in the Bott-Chern cohomology
group $H^2_{BC}(M,L)=H^2_{BC}(\tilde M)$. This implies that
$$V:=\mathrm{span}\{[\tom], \Ll_{A^c}[\tom]\}\subset H^2_{BC}(M,L)$$ is
$2$-dimensional. Then, obviously,
$\Ll_{A^c}$ acts on $V$ with two
1-dimensional eigenspaces, corresponding
to $\1\lambda$ and $-\1\lambda$. As $\Ll_{A^c}$ acts on $V$
essentially as a rotation with $\lambda\pi/2$, the flow of $A^c$, $e^{t
A^c}$,
will satisfy:
$$e^{t A^c}[\tilde \omega]=[\tilde \omega], \, \text{for}\,\, t=
2n\pi\lambda^{-1},\, n\in\ZZ.$$
We also note that
\begin{equation}\label{doi}
 \int_0^{2\pi\lambda^{-1}} e^{t A^c}[\tilde \omega]dt=0.
\end{equation}
Let now
\[
   \tilde \omega_W:= \int_0^{2\pi\lambda^{-1}} e^{t A^c}\tilde\omega dt.
\]
This new form is obtained as a sum of K\"ahler  forms
with the same automorphy, hence it is also an automorphic
K\"ahler  form. Its Bott-Chern class is equal to
$\int_0^{2\pi} e^{t A^c} [\tilde \omega] dt$, and thus
it vanishes by \eqref{doi}.

In conclusion, $\tom_W$ is a K\"ahler form
with trivial Bott-Chern class, and hence it admits a global automorphic
potential. We proved \ref{_S^1_main_Theorem_}.

\hfill

\remark
Another way to arrive at a K\"ahler form with potential is by averaging
using a kind of convolution. Let
$$
\psi= \begin{cases} \cos t+1, \, \,\, \,\, \text{for}\,\, t\in [-\pi, \pi]\\
0,  \qquad \qquad \text{for}\,\, t\not\in [-\pi, \pi].
\end{cases}
$$
Define
$$\tom_\psi=\int_\RR e^{t\lambda^{-1}A^c}\tom\psi(t)dt.$$
One can see that
$\Ll_{\lambda^{-1}A^c}\tom_\psi= \tom_{\psi'}$ and
$\Ll_{\lambda^{-1}A^c}^2\tom_\psi=\tom_{\psi''}$. Then, \eqref{unu} becomes
\begin{equation*}
 \begin{split}
  dd^c|A|_\psi^2&= \lambda^2\tilde \omega_\psi  +\Ll_{A^c}^2\tilde
\omega_{\psi}\\
&=\lambda^2(\tom_\psi+\tom_{\psi''})=
\lambda^2\int_\RR e^{t\lambda^{-1}A^c}\tom(\psi+\psi'')(t)dt.
 \end{split}
\end{equation*}
where $|A|_\psi^2$ means square length of $A$ taken
with respect to the metric $\omega_\psi$.
As $\psi''+\psi =1$ on $[-\pi, \pi]$, we see that  $ dd^c|A|_\psi^2>0$ and
hence $|A|_\psi^2$ is a K\"ahler potential for $\tom_\psi$. On the other
hand, one can verify that
\[
  dd^c(|A|^2_\psi)=\Ll_{A^c}^2\tilde\omega_\psi +
  \la^2\tilde\omega_{\psi''}= \tilde\omega_W,
\]
where $\tilde\omega_W=\int_{-\pi}^{\pi} e^{t A^c} \tilde \omega dt$.

Therefore, this averaging construction with ``weight'' $\psi$
gives the same form
$\tilde\omega_W=\int_{-\pi}^{\pi} e^{t A^c} \tilde \omega dt$
which we have obtained by the means of
averaging with the circle.

{\small

\noindent {\sc Liviu Ornea\\
University of Bucharest, Faculty of Mathematics, \\14
Academiei str., 70109 Bucharest, Romania. \emph{and}\\
Institute of Mathematics ``Simion Stoilow" of the Romanian Academy,\\
21, Calea Grivitei Street
010702-Bucharest, Romania }\\
\tt Liviu.Ornea@imar.ro, \ \ lornea@gta.math.unibuc.ro

\hfill

\noindent {\sc Misha Verbitsky\\
{\sc  Institute of Theoretical and
Experimental Physics \\
B. Cheremushkinskaya, 25, Moscow, 117259, Russia }\\
\tt verbit@maths.gla.ac.uk, \ \  verbit@mccme.ru
}
}

\end{document}